\documentclass[12pt,english]{amsart}
\usepackage[T1]{fontenc}
\pdfoutput=1
\usepackage[latin9]{inputenc}
\usepackage{graphicx}

\makeatletter

\providecommand{\tabularnewline}{\\}

\usepackage{amsmath}
\usepackage{amssymb}
\usepackage{amsthm}
\usepackage{latexsym}
\usepackage{diagrams}

\usepackage{lscape}

\usepackage{ulem}

\usepackage{pifont}
 \usepackage{hyperref}

\newtheorem{theorem}{Theorem}[section]
\newtheorem{corollary}[theorem]{Corollary}
\newtheorem{conjecture}[theorem]{Conjecture}

\newtheorem{definition}[theorem]{Definition}

\theoremstyle{remark}

\def\bc{\begin{center}}
\def\ec{\end{center}}
\def\be{\begin{enumerate}}
\def\ee{\end{enumerate}}
\def\bi{\begin{itemize}}
\def\ei{\end{itemize}}
\def\bs{\begin{slide}}
\def\es{\end{slide}}

\def\l[{\left[}
\def\r]{\right]}
\def\l({\left(}
\def\r){\right)}
\def\Aut{\text{Aut}}

\newcommand{\mc}[1]{\ensuremath{\mathcal{#1}}}
\newcommand{\Core}{\text{Core}}

\title{Representing the Sporadic Archimedean Polyhedra as Abstract Polytopes}

\author{Michael I. Hartley}
\address[Michael I. Hartley]{DownUnder GeoSolutions\\
80 Churchill Ave\\
Subiaco, 6008\\
Western Australia}
\email{mikeh@dugeo.com}
\author{Gordon I. Williams}
\address[Gordon I. Williams]{Department of Mathematics and Statistics\\ University of Alaska Fairbanks\\ PO Box 756660\\ Fairbanks, Alaska 99775-6660}
\urladdr{http://www.supposenot.com}
\email{giwilliams@alaska.edu}
\thanks{The authors would like to acknowledge the support of the Banff International Research Station and the organizers of the 2-day workshop on Convex and Abstract Polytopes: Ted Bisztriczky, Egon Schulte and Asia Ivic Weiss. We thank them for introducing us and for providing the time during which the seeds of this project were planted.}

\usepackage{babel}
\makeatother

\begin{document}
\normalem 

\begin{abstract}
We present the results of an investigation into the representations
of Archimedean polyhedra (those polyhedra containing only one type
of vertex figure) as quotients of regular abstract polytopes. Two
methods of generating these presentations are discussed, one of which
may be applied in a general setting, and another which makes use of
a regular polytope with the same automorphism group as the desired
quotient. Representations of the 14 sporadic Archimedean polyhedra
(including the pseudorhombicuboctahedron) as quotients of regular
abstract polyhedra are obtained, and summarised in a table. The information
is used to characterise which of these polyhedra have acoptic Petrie
schemes (that is, have well-defined Petrie duals). 
\end{abstract}
\maketitle

\section{Introduction}

Much of the focus in the study of abstract polytopes has been on the
study of the regular abstract polytopes. A publication of the first
author \cite{Har99} introduced a method for representing any abstract
polytope as a quotient of regular polytopes. In the current work
we present the application of this technique to the familiar, but
still interesting, Archimedean polyhedra and discuss implications
for the general theory of such representations that arose in trying
to systematically develop these representations. We discuss the theory
and presentations of the thirteen classical (uniform) Archimedean
polyhedra as well as the pseudorhombicuboctahedron, which we will
refer to as the fourteen sporadic Archimedean polyhedra.  In a separate
study, we will present and discuss the presentations for the two infinite
families of uniform convex polyhedra, the prisms and antiprisms.

\subsection{Outline of topics}

Section \ref{sec:APQP} reviews the structure of abstract polytopes
and their representation as quotients of regular polytopes and discusses
two new results on the structure of the quotient representations of
abstract polytopes. Section \ref{sec:RepFaces} describes a simple
method for developing a quotient presentation for a polyhedron from
a description of its faces. In Section \ref{sec:RepIso} we discuss
an alternative method of developing a quotient presentation for polytopes
that takes advantage of the structure of its automorphism group, and
in Section \ref{sec:IsoGeo} we develop this method more fully for
the specific polyhedra under study here. Finally, in Section \ref{sec:Analysis}
we discuss an example of how these quotient representations may be
used to answer questions about their structure computationally and
in Section \ref{sec:Conclusion} we present some of the open questions
inspired by the current work.

\section{Abstract Polytopes and Quotient Presentations}

\label{sec:APQP} To place the current work in the appropriate context
we must first review the structure of abstract polytopes and the central
results from the first author's \cite{Har99} for representing any
polytope as a quotient of regular abstract polytopes.

An {\em abstract polytope $P$ of rank $d$} (or {\em $d$-polytope})
is a graded poset with additional constraints chosen so as to generalize
combinatorial properties of the face lattice of a convex polytope.
Elements of these posets are referred to as {\em faces}, and a
face $F$ is said to be {\em contained} in a face $G$ if $F<G$
in the poset. One consequence of this historical connection to convex
polytopes is that contrary to the usual convention for graded posets,
the rank function $\rho$ maps $P$ to the set $\{-1,0,1,2,...,d\}$
so that the minimal face has rank $-1$, but otherwise satisfies the
usual conditions of a rank function. A face at rank $i$ is an {\em
$i$-face}. A face $F$ is {\em incident} to a face $G$ if either
$F<G$ or $G<F$. A {\em proper face} is any face which is not
a maximal or minimal face of the poset. A {\em flag} is any maximal
chain in the poset, and the {\em length} of a chain $C$ we define
to be $|C|-1$. Following \cite{McMSch02} we will require that the
poset $P$ also possess the following four properties:

\begin{description}
\item [{P1}] $P$ contains a least face and a greatest face, denoted $F_{-1}$
and $F_{d}$ respectively; 
\item [{P2}] Every flag of $P$ is of length $d+1$; 
\item [{P3}] $P$ is strongly connected; 
\item [{P4}] For each $i=0,1,...,d-1$, if $F$ and $G$ are incident faces
of $P$, and the ranks of $F$ and $G$ are $i-1$ and $i+1$ respectively,
then there exist precisely two $i$-faces $H$ of $P$ such that $F<H<G$. 
\end{description}
Note that an abstract polytope is {\em connected} if either $d\le1$,
or $d\ge2$ and for any two proper faces $F$ and $G$ of $P$ there
exists a finite sequence of incident proper faces $J_{0},J_{1},...,J_{m}$
such that $F=J_{0}$ and $G=J_{m}$. A polytope is {\em strongly
connected} if every section of the polytope is connected, where a
{\em section } corresponding to the faces $H$ and $K$ is the
set $H/K:=\{F\in P\mid H<F<K\}$. Some texts are more concerned with
the notion of {\em flag connectivity}. Two flags are {\em adjacent}
if they differ by only a single face. A poset is {\em flag-connected}
if for each pair of flags there exists a sequence of adjacent flags
connecting them, and a poset is {\em strongly flag-connected} if
this property holds for every section of the poset. It has been shown
\cite{McMSch02} that for any poset with properties \textbf{P1} and
\textbf{P2}, being strongly connected is equivalent to being strongly
flag-connected. A polytope is said to be {\em regular} if its automorphism
group $\mathrm{Aut({\mc P})}$ acts transitively on the set ${\mc F}({\mc P})$
of its flags.

To understand what follows, a basic understanding of the structure
of string C-groups is necessary, so we will review the essential definitions
here. A {\em C-group} $W$ is a group generated by a set of (distinct)
involutions $S=\{s_{0},s_{1},\ldots,s_{n-1}\}$ such that $\langle s_{i}|i\in I\rangle\cap\langle s_{j}|j\in J\rangle=\langle s_{i}|i\in I\cap J\rangle$ for all $I,J$
(the so-called {\em intersection property}). Coxeter groups are
the most famous examples of C-groups (see \cite{Hum90},\cite{McMSch02}).
A C-group is a {\em string} C-group if $(s_{i}s_{j})^{2}=1$ for
all $|i-j|>1$. An important result in the theory of abstract polytopes
is that the regular polytopes are in one-to-one correspondence with
the string C-groups, in particular, that the automorphism group of
any regular abstract polytope is a string C-group and that from every
string C-group $W$ a unique regular polytope ${\mc P}(W)$ may be
constructed whose automorphism group is $W$ \cite{McMSch02}. Given
a C-group $W$ and a polytope ${\mc Q}$ (not necessarily related
to ${\mc P}$), we may attempt to define an action of $W$
on ${\mc F}({\mc Q})$ as follows. For any flag $ \Phi$ of
${\mc Q}$, let $\Phi^{s_{i}}$ be the unique flag differing from
$\Phi$ only by the element at rank $i$. If this extends to a well-defined
action of $W$ on ${\mc F}({\mc Q})$, it is called the {\em flag
action} of $W$ on (the flags of) ${\mc Q}$. The flag action should
not be confused with the natural action of the automorphism group
$W$ of a regular polytope \mc{Q} on its flags. As noted in \cite{Har99},
it is always possible to find a C-group acting on a given abstract
polytope $\mc Q$ (regular or not) via the flag action.

We consider now the representation of abstract polytopes first presented
as Theorem 5.3 of \cite{Har99}. \begin{theorem}\label{th:QuotRep}
Let $\mathcal{Q}$ be an abstract $n$-polytope, $W$ any string C-group
acting on the flags of $\mathcal{Q}$ via the flag action and $\mathcal{P}(W)$
the regular polytope with automorphism group $W$. If we select any
flag $\Phi$ as the base flag of $\mathcal{Q}$ and let $N=\{a\in W\mid\Phi^{a}=\Phi\}$, then $\mathcal{Q}$ is isomorphic to $\mathcal{P}(W)/N$. Moreover,
two polytopes are isomorphic if and only if they are quotients $\mathcal{P}(W)/N$
and $\mathcal{P}(W)/N'$ where $N$ and $N'$ are conjugate subgroups
of $W$. \end{theorem}

An interesting fact about these presentations that does not seem to
appear explicitly elsewhere in the literature is that there is a strong
relationship between the number of transitivity classes of flags under the automorphism group in
the polytope and the number of conjugates of the stabilizer subgroup
$N$. This relationship is formalized as follows.
 \begin{theorem}The
number of transitivity classes of flags under the automorphism group in a polytope ${\mc Q}$
is equal to the number of conjugates in $W$ of the stabilizer subgroup
$N$ for any choice of base flag $\Phi$ in its quotient presentation,
that is, $|W:\mathrm{Norm}_{W}(N)|$. \label{th:flagCount} \end{theorem} 

\begin{proof}
Let $\Phi$ and $\Phi'$ be two flags of a polytope $\mc Q$, let
$W$ be a string C-group acting on \mc{Q}, and let \mc{P} be the
regular polytope whose automorphism group is $W$ (so $\mc P=\mc P(W)$).
Let $N$ be the stabilizer of $\Phi$ in $W$, and let $N'$ be the
stabilizer of $\Phi'$ in $W$. Let $\Phi'=\Phi^{u}$, so that $N'=N^{u}$.
Let $\psi$ be an automorphism of \mc{Q} with $\Phi\psi=\Phi'$,
and suppose $n\in N$. Observe then that \begin{align*}
(\Phi')^{n} & =(\Phi\psi)^{n} &  & \text{by the definition of $\psi$}\\
 & =(\Phi^{n})\psi &  & \text{by Lemma 4.1 of \cite{Har99}}\\
 & =(\Phi)\psi &  & \text{since $n\in N$, the stabilizer of $\Phi$}\\
 & =\Phi' &  & \text{by the definition of $\psi$.}\end{align*}
 Therefore, $n\in N'$, so $N=N'$.

Conversely, let $N=N'$. Then, a map from $\mc P/N$ to $\mc P/N'$
may be constructed as in the proof of Theorem 5.3 of \cite{Har99}
(our Theorem \ref{th:QuotRep}), which does indeed map $\Phi$ to $\Phi'$. 
\end{proof}
Theorem \ref{th:QuotRep} does not provide much guidance on finding
an efficient (i.e. small) presentation for a given polytope. In particular,
it is interesting to try to determine what the smallest regular polytope
is that may be used as a cover of a given polytope under the flag action of the automorphism group of the regular polytope. Let $\Core(W,N)$
be the subgroup of $N$ obtained as $\bigcap\limits _{w\in W}N^{w}$,
in other words, the largest normal subgroup of $W$ in $N$.

\begin{theorem} Let $\mc P(W/\Core(W,N))$ be a well defined regular
polytope, and $\mc R$ any other regular cover of $\mc P(W)/N$ whose automorphism group acts on $\mc P(W)/N$ via the flag action, and on which $W$ acts likewise.
Then $\mc R$ also covers $\mc P(W/\Core(W,N))$.\label{th:core}\end{theorem} 

\begin{proof}
Let $\mc R=\mc P(W)/K=\mc P(W/K)$ be a regular cover for $\mc P(W)/N$.
Then the flag action of $W/K$ on $\mc P(W)/N$ is well defined; that
is, for any $w\in W$ and any flag $\Phi$ of $\mc P(W)/N$, we have
$\Phi^{wK}$ is well defined, because $\Phi^{wk}$ independent of
the choice of $k$ in $K$, but depends only on $w$. It follows that
for all $k\in K$, any $w\in W$, and any flag $\Phi$ of $\mc P(W)/N$,
we have $(\Phi^{wk})^{w^{-1}}=\Phi$, so, $wkw^{-1}\in N$. Therefore,
$k\in N^{w}$ for all $w\in W$, so $k\in\Core(W,N)$. 
\end{proof}
Now, Theorem 3.4 of \cite{Har99a} states that \[
\Gamma(\mc P(W)/N)\cong W/\Core(W,N),\]
 where $\Gamma(\mc P(W)/N)$ is the image of the homomorphism induced
by the flag action from $W$ into $Sym(Flags(\mc P(W)/N))$. In the
case that $\mc P(W)/N$ is a finite polytope, so that $N$ has finite
index in $W$, it follows that $\Core(W,N)$ is a finite index
normal subgroup of $W$. This is because $W$ acts on the finitely
many right cosets of $N$ via right multiplication, leading to a homomorphism
from $W$ to $\Sigma=Sym(|W:N|)$. The kernel of this homomorphism
is $\Core(W,N)$, and thus $W/\Core(W,N)$ is isomorphic to a subgroup
of the finite group $\Sigma$. Hence, a finite polytope always has
a finite regular cover if $W/\Core(W,N)$ is a C-group. No proof that $W/\Core(W,N)$ is indeed a C-group has yet been published.

Barry Monson notes (\cite{Mon}) that there exist quotients $\mc{Q}=\mc{P}/N$ of a polytope \mc{P}, for which the flag action of the automorphism group $W$ of \mc{P} on \mc{Q} is not well defined. The theory of such exceptional quotients is not well developed. This article therefore concerns itself exclusively with quotients of \mc{P} on which the flag action of $\Aut(\mc{P})$ is well-defined.

\section{An Example in Detail}

\label{sec:RepFaces} 
\begin{figure}[htbp]

\begin{centering}
\includegraphics[width=3in]{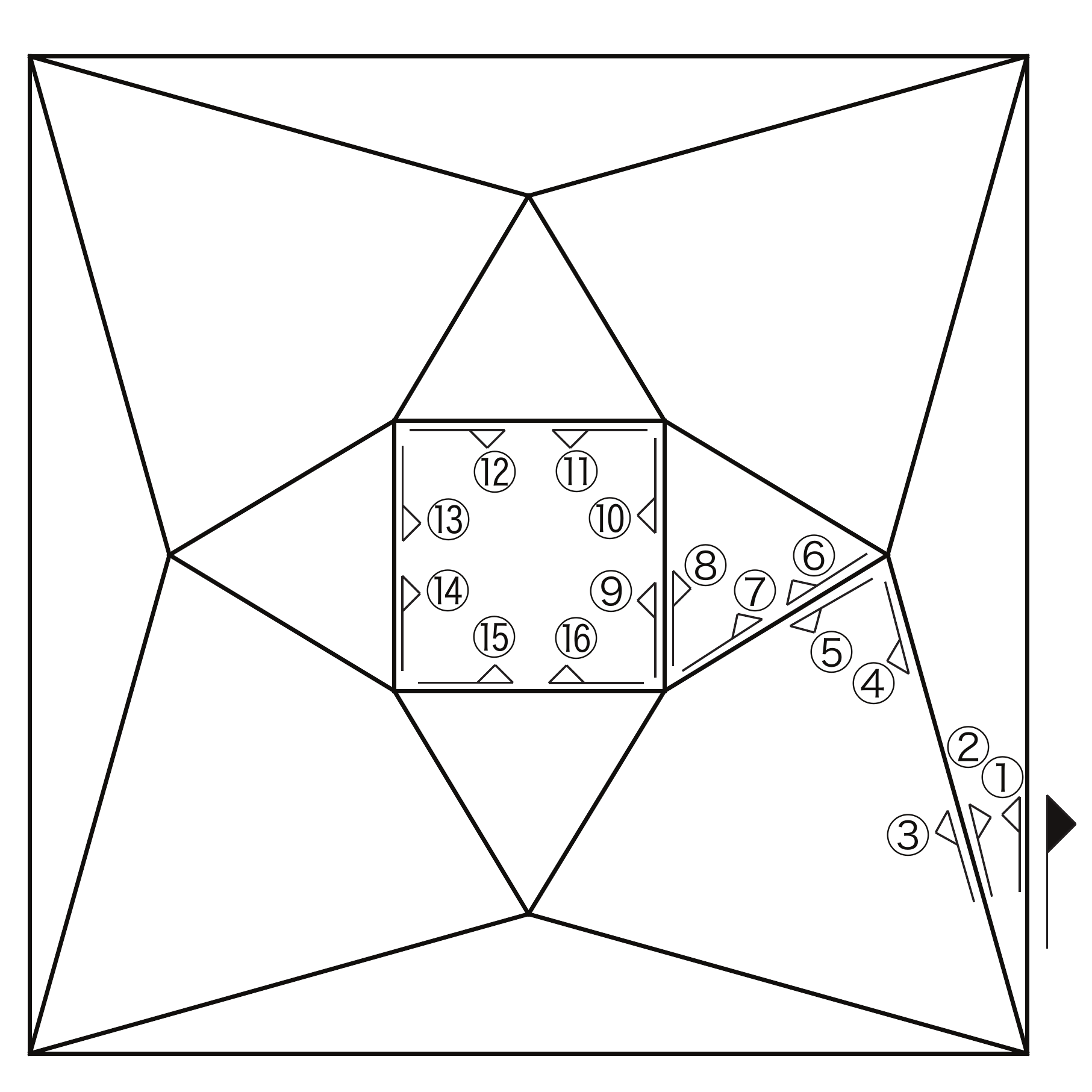}\includegraphics[width=3in]{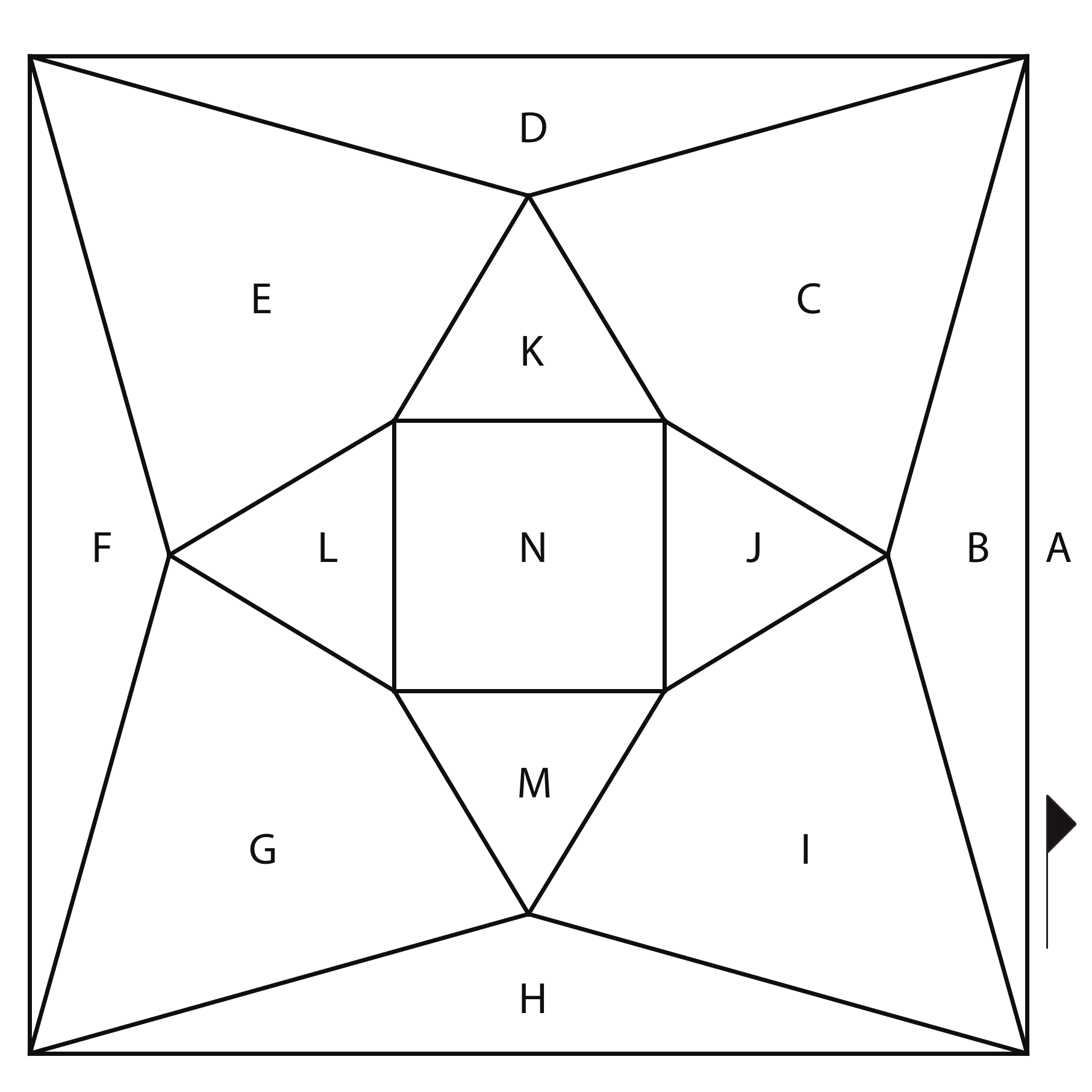} 
\par\end{centering}

\caption{On the left is pictured the cuboctahedron with a sequence of labeled
flags used in the construction of the stabilizer subgroup of the base
flag indicated in black. On the right is the same figure with labels
indicated for each of the faces of the cuboctahedron.}

\label{fig:cuboct} 
\end{figure}

From Theorem \ref{th:QuotRep} we learn that
any given polytope \mc{Q} admits a presentation as the quotient
of a regular polytope. To find such a presentation we must first identify
a string C-group $W$ acting on the flags of \mc{Q} via the flag
action, and then having selected a base flag $\Phi\in\mc Q$, we must
identify the stabilizer of $\Phi$ in $W$. To illustrate the mechanics
of this process we will consider here the case of the cuboctahedron.
As in \cite{Gr03} we will associate to each uniform or Archimedean
polyhedron a symbol of type $p_{1}.p_{2}...p_{k}$, which specifies
an oriented cyclic sequence of the number of sides of the faces surrounding
each vertex. For example, 3.4.3.4 designates the cuboctahedron, which
is an isogonal polyhedron with a triangle, a square, a triangle and
a square about each vertex in that cyclic order. Figure \ref{fig:cuboct}
shows the corresponding graph of the one-skeleton of the cuboctahedron.

First we select as our C-group the group $W=\langle a,b,c\mid a^{2}=b^{2}=c^{2}=(ac)^{2}=(ab)^{12}=(bc)^{4}=e\rangle$,
where $e$ is the identity. For ease of notation we write $a,b,c$ instead of $s_{0},s_{1},s_{2}$, respectively. In general, one possible choice of the string C-group acting
on a 3-polytope is the group ${W=\langle a,b,c\mid a^{2}=b^{2}=c^{2}=(ac)^{2}=(ab)^{j}=(bc)^{k}=e\rangle}$,
where $j$ is the least common multiple of the number of sides of
polygons in \mc{Q} and $k$ is the least common multiple of the
degrees of the vertices of \mc{Q}. Here the action of the generators
$a,b$ or $c$ on a flag $F$ of \mc{Q} yields the adjacent flag
differing from $F$ only by the vertex, edge or face, respectively.

For our choice of base flag in this example we select a flag $\Phi$
on a square face (in our diagram this corresponds to the outside face),
and we mark it with a solid black flag. Construction of the stabilizer
subgroup $N$ of $\Phi$ in $W$ is a bit more involved. For each
of the faces of \mc{Q} we may construct a sequence of consecutively adjacent flags
starting at the base flag, going out to the face, forming a circuit
of the edges and vertices of the face, and returning to the base flag.
Each of these flags may be obtained from $\Phi$ via the flag action of $W$ on
$\Phi$; for example, the flag marked with a \ding{172} is obtained
from $\Phi$ via the action of the generator $c$ of $W$. Starting
at flag \ding{180}, a complete circuit of the face \textsf{N} is
obtained from flag \ding{180} by application of the element $(ab)^{4}\in W$.
Thus the group element corresponding to starting at the base flag
and traversing the face marked \textsf{N} and returning is $((ab)^{4})^{cbacbacbc}$.

Let $N$ be the group in $W$ generated by \begin{multline}
\{(ab)^{4},((ab)^{3})^{c},((ab)^{4})^{cbabc},((ab)^{3})^{cba},((ab)^{4})^{cbcabab},\\
((ab)^{3})^{cbab},((ab)^{4})^{cbacb},((ab)^{3})^{cb},((ab)^{4})^{cbc},((ab)^{3})^{cbcabc},\\
((ab)^{3})^{cbcabcba},((ab)^{3})^{cbcabcabab},((ab)^{3})^{cbabacbc},((ab)^{4})^{cbacbacbc}\}\label{eq:generators}\end{multline}
 The generators in  (\ref{eq:generators}) correspond to
faces \textsf{A} through \textsf{N} in Figure \ref{fig:cuboct} in
that order.

Note that in general, finding elements of $W$ that, as above, traverse each
face of $\mc Q$ may only suffice to generate a proper subgroup of
$N$. Inspection of $\mc Q$ should then reveal other elements of
$W$ that stabilise $\Phi$~-- these can then be added to the generating
set for $N$. In the example here, however, the elements listed do
indeed generate the whole of the base flag stabiliser $N$.  Then by Theorem \ref{th:QuotRep} the cuboctahedron \mc{Q}
is isomorphic to $\mc P(W)/N$.

\section{Representation via Isomorphism}

\label{sec:RepIso} In the context of the current work, an important
observation is that the automorphism group of a polyhedron is often
shared with a better understood regular polytope. For example, the
automorphism group of the cuboctahedron is that of the cube. It turns
out that the quotient presentation can be characterized with the help  of
the symmetry group of the associated regular polytope. Again, we let
$\mc P$ be a regular $n$-polytope, with automorphism group $W$. Let
$\mc Q$ be a quotient $\mc P/N$ of $\mc P$ (not necessarily regular) admitting the flag action by $W$ with $\Psi$ a base flag for \mc Q chosen so that $N$ is the stabilizer for $\Psi$,
and let $\mc R$ be a regular $d$-polytope whose automorphism group is
isomorphic to $\Aut(Q)$. Note that we do not assume that $d=n$. Let $\Aut(\mc R)=\langle\rho_{0},\rho_{1},\dots,\rho_{d-1}\rangle$.
Let $\phi$ be an isomorphism from $\Aut(\mc R)$ to $\Aut(\mc P/N)$.

Let $\Phi$ be a flag of $\mc R$, and $\Psi$ a base flag of $\mc Q=\mc P/N$,
stabilized by $N$ under the flag action. For each $\rho_{i}$, let
$\nu_{i}$ be an element of $W$ that maps $\Psi$ to $\Psi(\rho_{i}\phi)$
under the flag action, that is, $\Psi^{\nu_{i}}=\Psi(\rho_{i}\phi)$.
Let $V$ be the subgroup of $W$ generated by the $\nu_{i}$. Finally,
define a map $\psi$ taking words $w$ in the generators of $\Aut(\mc R)$
to the group $W$, via $w\psi=(\rho_{i_{1}}\rho_{i_{2}}\dots\rho_{i_{k}})\psi=\nu_{i_{k}}\dots\nu_{i_{2}}\nu_{i_{1}}$. Note that the action of $\psi$ reverses the order of the generators.

The following result goes a long way towards characterizing $N$ in
terms of $\Aut(\mc R)$.

\begin{theorem} The set $N\cap V$ is the set of all images $w\psi$
of words $w$ in the $\rho_{i}$ such that $w=1$ as an element of
$\Aut(\mc R)$. \label{th:w=00003D1}\end{theorem} 

\begin{proof}
Note that $\rho_{i_{1}}\dots\rho_{i_{k}}=1$ in $\Aut(\mc R)$ if
and only if $\Psi((\rho_{i_{1}}\dots\rho_{i_{k}})\phi)=\Psi$. This
will be so if and only if $\Psi(\rho_{i_{1}}\phi)\dots(\rho_{i_{k}}\phi)=\Psi$.
Since the flag action commutes with the action of the automorphism
group (Lemma 4.1 of \cite{Har99}), we have \begin{align*}
\left(\Psi(\rho_{i_{j}}\phi)\dots(\rho_{i_{k}}\phi)\right)^{\nu_{i_{j-1}}\dots\nu_{i_{1}}} & =\left(\Psi^{\nu_{i_{j}}}(\rho_{i_{j+1}}\phi)\dots(\rho_{i_{k}}\phi)\right)^{\nu_{i_{j-1}}\dots\nu_{i_{1}}}\\
 & =\left(\Psi(\rho_{i_{j+1}}\phi)\dots(\rho_{i_{k}}\phi)\right)^{\nu_{i_{j}}\dots\nu_{i_{1}}}.\end{align*}
 Thus, $\Psi(\rho_{i_{1}}\phi)\dots(\rho_{i_{k}}\phi)=\Psi$ if and
only if $\Psi^{\nu_{i_{k}}\dots\nu_{i_{1}}}=\Psi$, that is, if and
only if $\nu_{i_{k}}\dots\nu_{i_{1}}=w\psi\in N$. This completes
the proof.
\end{proof}
So the elements of $N\cap V$ have been characterized. To characterize
the whole of $N$, it is sufficient to characterize elements of $N\cap V\mu$,
for arbitrary cosets $V\mu$ of $V$ in $W$. This is not as difficult
as it may seem. Note that if $\mu\in N$, then $N\cap V\mu=(N\cap V)\mu$.

\begin{theorem} Let $T$ be a right transversal of $V$ in $W$,
such that for all $\mu\in T$, if $N\cap V\mu\neq\emptyset$, then
$\mu\in N$. Then \[
N=\bigcup\limits _{\mu\in N\cap T}\left\{ (w\psi)\mu:w=1{\mathrm{~in~}}\Aut(\mc R)\right\} .\]
\label{th:N=00003Dcup} \end{theorem} 

\begin{proof}
For any right transversal $T$ of $V$ in $W$, \[
N=N\cap W=N\cap\left(\bigcup_{\mu\in T}V\mu\right)=\bigcup_{\mu\in T}\left(N\cap V\mu\right).\]
 For the transversal chosen here, $N\cap V\mu$ is empty unless $\mu\in N$,
whence also $N\cap V\mu=(N\cap V)\mu$. It follows that \[
N=\bigcup_{\mu\in N\cap T}\left((N\cap V)\mu\right),\]
 which by Theorem \ref{th:w=00003D1} is \[
N=\bigcup_{\mu\in N\cap T}\left\{ (w\psi)\mu:w=1{\mathrm{~in~}}\Aut(\mc R)\right\} \]
 as desired. 
\end{proof}
This gives a characterisation of the elements of $N$, in terms of
the elements of $\Aut(\mc R)$, the map $\phi$, and the transversal
$T$.

Theorems \ref{th:w=00003D1} and \ref{th:N=00003Dcup} are particularly
useful for the purposes of this article, since every uniform sporadic
Archimedean solid has an automorphism group that is also the automorphism
group of a regular polytope $\mc R$. In most cases, the choice of
$\mc R$ is obvious~--- it will be the underlying platonic solid.
The snub cube and snub dodecahedron have as automorphism group the
\textit{rotation} group of the cube and dodecahedron respectively,
not the full automorphism groups. However, these rotation groups are
isomorphic (respectively) to the automorphism groups of the hemi-cube
$\{4,3\}_{3}$ and the hemi-dodecahedron $\{5,3\}_{5}$, so these
theorems may still be applied.

In the following sections, Theorems \ref{th:w=00003D1} and \ref{th:N=00003Dcup}
are used to construct each of the Archimedean solids as a quotient
$\mc P/N$ of some regular polytope $\mc P$ by a subgroup $N$ of
its automorphism group. The steps in construction are as follows.

\begin{enumerate}
\item Find a polytope $\mc P$ that is known to cover the desired Archimedean
solid. 
\item Identify, using Theorems \ref{th:w=00003D1} and \ref{th:N=00003Dcup},
a subset $S$ of $N$. 
\item Prove, or computationally verify, that $S$ generates a subgroup of
$\Aut(\mc P)$ whose index is the same as the (known) index of $N$. 
\item Finally, use Theorem~\ref{th:core} to find a minimal regular cover
$\mc P/\Core(\Aut(\mc P),N)$ for the Archimedean solid $\mc P/N$. 
\end{enumerate}
The index of $N$ in $\Aut(\mc P)$ is known, from Theorem 2.5 of
\cite{Har99a}, to be just the number of flags of the quotient $\mc P/N$,
which is easy to compute. Indeed, the Archimedean solid with symbol $p_{1}.p_{2}\ldots p_{k}$ has exactly $2k$ flags at every vertex.

\section{Isomorphism for Geometric Operations}

\label{sec:IsoGeo}

\begin{figure}
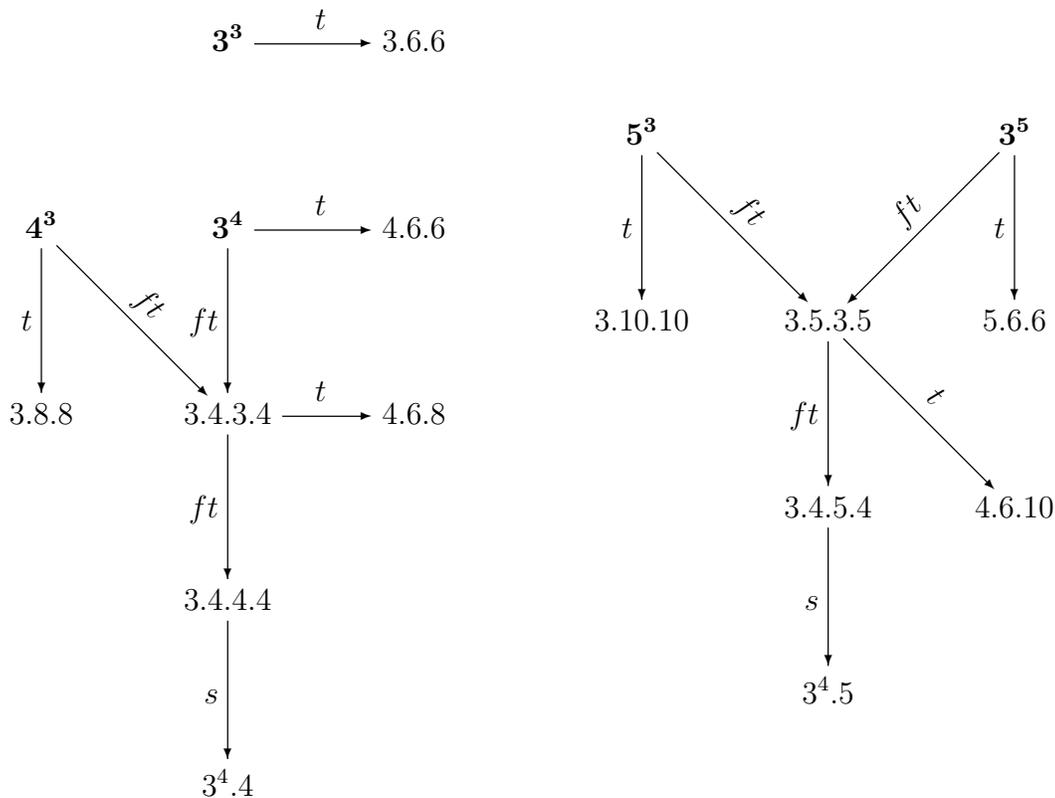

$\begin{array}{ccc}
\begin{diagram}
		 & 		& \mathbf{3^{3} }		&\rTo^{t}	&3.6.6\\
 		& 		&    	&		&\\
 \mathbf{4^{3}}	& 		&  \mathbf{3^{4} }		&\rTo^{t}	&4.6.6\\
 \dTo^{t}	&\rdTo^{ft}	&\dTo^{ft}		&		&\\
 3.8.8	&		&3.4.3.4		&\rTo^{t}	&4.6.8\\
 		&		&\dTo^{ft}		&		&\\
		&		&3.4.4.4		&		&\\
		&		&\dTo^{s}		&		&\\
		&		&3^{4}.4		&		&
\end{diagram}
&\hspace{.5in} &
\begin{diagram}
\mathbf{5^{3}	}	&		&		&		&\mathbf{3^{5}}\\
\dTo^{t}	&\rdTo^{ft}	&		&\ldTo^{ft}	&\dTo^{t}\\
3.10.10	&		&3.5.3.5	&		&5.6.6\\
		&		&\dTo^{ft}	&	\rdTo^{{t}}	&\\
		&		&3.4.5.4	&		&4.6.10\\
		&		&\dTo^{s}	&		&\\
		&		&3^{4}.5	&		&
\end{diagram}\end{array}$
\caption{The construction of the sporadic uniform Archimedean polyhedra from the Platonic solids.}\label{f:inheritance}
\end{figure}

From a combinatorial -- but not geometric -- standpoint, each of the
uniform sporadic Archimedean polyhedra may be constructed from a
Platonic solid by (possibly repeated) application of either truncation,
full truncation, rhombification or snubbing. By Theorem \ref{th:w=00003D1},
we may construct quotient presentations for these polyhedra by determining
the appropriate choices for the $\nu_{i}\in V$ that correspond to
these operations. Let us now carefully define what each of these operations
does. Geometrically, {\em truncation} ($t$) cuts off each of the
vertices of the polyhedron, replacing them with the corresponding
vertex figure as a facet. {\em Full truncation} ($ft$) performs
essentially the same operation, but the cut is taken deeper so that
new facets share a vertex if the corresponding vertices shared an
edge, and all of the original edges are replaced with single vertices.
{\em Rhombification} ($r$) is a little more complicated geometrically,
but from a combinatorial standpoint is equivalent to applying full
truncation twice (the difficulty is in getting the new facets to be
geometrically regular). Finally, to construct the {\em snub} of
a polyhedron requires first constructing the rhombification, and then
triangulating the squares generated by the second full truncation
in such a way as to preserve the rotational symmetries of the figure
(in Figure \ref{f:inheritance} the triangulation step is indicated
by $s$). The ways in which each of the sporadic uniform Archimedean
polyhedra may be obtained (hierarchically) from the Platonic solids
via these operations is given in Figure \ref{f:inheritance}. Note for instance that $4^{3}$ abbreviates the symbol $4.4.4$ for the cube. More information on these, and other, operations on the maps associated with polyhedra is available in \cite{PisRan00}.

\subsection{Generators of $V$}

For the convenience of the reader, we present here the morphisms $\psi$ 
from the words in the generators of the symmetry groups of the regular polyhedra $\mc R$ into the symmetry
groups of the regular covers $\mc P$ of the quotient polytopes $\mc Q$ that provide
the generators for the subgroup $V$ of Theorems \ref{th:w=00003D1}
and \ref{th:N=00003Dcup}. It is also important to note that different
morphisms (and corresponding sets of generators) arise if one
makes different choices for the base flag in the quotient polytope
than those made here, and that $\nu_{0}$ and $\nu_{2}$ may be interchanged
by using the dual choice for the polytope $\mc R$ (where possible
and appropriate). The map $\psi$ in each case is determined by its
action on the generators of $\Aut(\mc R)$, denoted $\rho_{0},\rho_{1}$
and $\rho_{2}$, in terms of the generators of $W={\mc P}=\langle{a,b,c}\rangle$ in the usual way.

\subsubsection{Truncation}

There are five Archimedean polyhedra obtained by truncation of each
of the Platonic solids, namely, the truncated tetrahedron, cube, octahedron,
icosahedron and dodecahedron. In each instance the vertex star contains
either two hexagons, two octagons or two decagons. Here we choose
as a base flag $\Psi$ on one of those hexagons, octagons or decagons
whose edge is shared with another polygon of the same type. Thus \begin{align*}
\rho_{0}\psi=\nu_{0} & =a,\\
\rho_{1}\psi=\nu_{1} & =bab\\
\rho_{2}\psi=\nu_{2} & =c,\end{align*}
 so $V=\langle a,bab,c\rangle$.

\subsubsection{Full Truncation}

Full truncation provides derivations for two of the Archimedean polyhedra,
the cuboctahedron and the icosidodecahedron. We have chosen to perform
full truncation to the cube and the dodecahedron, respectively, and our base flags
on square or pentagonal faces respectively. Thus \begin{align*}
\rho_{0}\psi=\nu_{0} & =b,\\
\rho_{1}\psi=\nu_{1} & =a,\\
\rho_{2}\psi=\nu_{2} & =cbc,\end{align*}
 so $V=\langle b,a,cbc\rangle$.

\subsubsection{Rhombification}

There are two Archimedean polyhedra obtained by rhombification, the
small rhombicuboctahedron and the small rhombicosidodecahedron. Here
we begin with the cube and dodecahedron, respectively, and our base flag lies on
an edge of a square or pentagonal face shared with the square face
introduced by the second full truncation. Thus \begin{align*}
\rho_{0}\psi=\nu_{0} & =a,\\
\rho_{1}\psi=\nu_{1} & =b,\\
\rho_{2}\psi=\nu_{2} & =cbabc,\end{align*}
 so $V=\langle a,b,cbabc\rangle$. While
it is true that the octahedron may be obtained by full truncation
from the tetrahedron (and so the cuboctahedron may be obtained by
rhombification of the tetrahedron), the maps given do not provide
an isomorphism since the symmetry group of the octahedron, and hence
the cuboctahedron, is larger than that of the tetrahedron.

\subsubsection{Truncation of Full Truncation}

There are two Archimedean polyhedra obtained in this way, the great
rhombicuboctahedron and the great rhombicosidodecahedron. Here we
begin with a cube and a dodecahedron, respectively, and our base flag lies on either an octagonal or decagonal face with an edge shared with a square. Thus \begin{align*}
\rho_{0}\psi=\nu_{0} & =a,\\
\rho_{1}\psi=\nu_{1} & =bab,\\
\rho_{2}\psi=\nu_{2} & =cbabc\end{align*}
 and so $V=\langle a,bab,cbabc\rangle$

\subsubsection{Snubbing}

There are two Archimedean polyhedra obtained by the snubbing operation,
the snub cube and the snub dodecahedron. For the presentation given
below for $V$, we have chosen to start with the hemi-cube and the
hemi-dodecahedron, respectively. These regular polyhedra are non-orientable, so
the group of $\mc R$ is coincides with its rotation subgroup, and we need
only consider the generators of this group in determining $V$. In
each case the base flag lies on either a square or pentagonal face.
\begin{align*}
\rho_{1}\rho_{0}\psi=\nu_{0}\nu_{1} & =ab,\\
\rho_{2}\rho_{1}\psi=\nu_{1}\nu_{2} & =bcbabcbc\end{align*}
 so $V=\langle ab,bcbabcbc\rangle$.

\subsection{The Cuboctahedron}

To better understand how this works in practice, let us return to
the example of the cuboctahedron, conceived as the full truncation of
the cube. In this case $\Aut(\mc R)=\langle s,t,u\mid s^{2}=t^{2}=u^{2}=(su)^{2}=(st)^{4}=(tu)^{3}\rangle$,
and $V=\langle b,a,cbc\rangle<W$ (this $W$ was defined in Section
\ref{sec:RepFaces}). By Theorem \ref{th:w=00003D1} (and Theorem \ref{th:N=00003Dcup} if necessary), if we can find
a set of words in the generators $s,t,u$ of $\Aut(\mc R)$ that are equivalent to the identity
in $\Aut(\mc R)$ and whose images generate a group of the appropriate
index (in this case 96) in $W$, then we will have found the necessary
subgroup of $W$ for use in the quotient presentation of the cuboctahedron.
Recall that if $\phi$ is the isomorphism from $\Aut(\mc R)$ to $\Aut(\mc P/N)$, and $\psi$ the associated map from $\Aut(\mc R)$ to $W$, 
then $s\psi=b,t\psi=a$ and $u\psi=cbc$; using this map we
generate the list of words given below in Equation \ref{eq:newgen},
which satisfies the conditions of Theorem \ref{th:w=00003D1}: \begin{multline}
\{(st)^{4},(ut)^{3},((st)^{4})^{utu},((ut)^{3})^{st},((st)^{4})^{utsts},\\
((ut)^{3})^{sts},((st)^{4})^{uts},((ut)^{3})^{s},((st)^{4})^{u},((ut)^{3})^{stu},\\
((ut)^{3})^{stus},((ut)^{3})^{stutsts},((ut)^{3})^{ststu},((st)^{4})^{utstu}\}.\label{eq:newgen}\end{multline}
 Each of the terms in Equation \ref{eq:newgen} corresponds to either
a circuit of one of the square faces of the cube, or to a traversal
of one of the vertex stars of the cube (starting at, and returning
to a chosen base flag), and so clearly is equivalent to 1 in $\Aut(\mc R)$. By Theorem \ref{th:w=00003D1}, if we apply $\psi$ to each of these terms we obtain an element of the subgroup $N$ required to construct a quotient representation under the flag action of $W$. Conveniently, in this example
each of the terms in Equation \ref{eq:newgen} corresponds to one
of the generators given in Equation \ref{eq:generators} and are listed
in the same order. To see this, consider for example the sixth item
on the list, $((ut)^{3})^{sts}$. When we apply the map $\psi$, we
see that \begin{align*}
((ut)^{3})^{sts}\psi & =(stsutututsts)\psi\\
 & =babacbcacbcacbcbab&  & \text{(by definition of $\psi$)}\\
 & =bab{\bf ca}b{\bf ac}cb{\bf ac}cbcbab &  & \text{(by commutivity of $a$ and $c$ in $W$)}\\
 & =babcabababcbab=((ab)^{3})^{cbab} &  & \text{(since $c^{2}=1$)}\end{align*}
 as was desired.

We conclude this discussion with the results of constructing such presentations for each of the sporadic uniform Archimedean solids.
\begin{theorem} Each of the sporadic uniform Archimedean solids has a finite regular cover whose automorphism group acts on the Archimedean solid via the flag action. Moreover, the regular covers are minimal in this sense, as detailed in Table \ref{T:summary}.\end{theorem}
%
\begin{table}[t]

\caption{This summarizes the representations of the Archimedean solids as quotients
of abstract regular polytopes ${\mc P}={\mc P}(W)$. These $\mc P$ are the minimal regular polytopes whose automorphism groups act on the Archimedean solids via the flag action.}

\begin{centering}
\begin{tabular}{|c|c|c|c|c|c|}
\hline 
\textbf{Polytope}  & \textbf{Vertex}  & \textbf{Schläfli}  & $|W|$ & $|N|$ \tabularnewline
 & \textbf{Figure}  & \textbf{Type of ${\mc P}(W)$} &  &  \tabularnewline
\hline
\hline 
Trunc. Tetrahedron & 3.6.6  & $\{6,3\}$  & 144 & 2 \tabularnewline
\hline 
Trunc. Octahedron  & 4.6.6  & $\{8,3\}$  & 6912 & 48 \tabularnewline
\hline 
Cuboctahedron  & 3.4.3.4  & $\{12,4\}$  & 2304  & 24 \tabularnewline
\hline 
Trunc. Cube  & 3.8.8  & $\{24,3\}$  & 82944  & 576  \tabularnewline
\hline 
Icosadodecahedron  & 3.5.3.5  & $\{15,4\}$  & 14400  & 120 \tabularnewline
\hline 
Trunc. Icosahedron  & 5.6.6  & $\{30,3\}$  & 2592000  & 7200 \tabularnewline
\hline 
Sm. Rhombicuboctahedron  & 3.4.4.4  & $\{12,4\}$  & 1327104 & 6912  \tabularnewline
\hline
Pseudorhombicuboctahedron& 3.4.4.4 & $\{12,4\}$ &$ 2^{35}3^{5}5^{2}7\cdot 11$ &$2^{29}3^{4}5^{2}7\cdot 11$\tabularnewline
\hline
Snub Cube  & 3.3.3.3.4  & $\{12,5\}$  & $2^{32}3^{11}5^{1}$ & $2^{28}3^{10}$ \tabularnewline
\hline 
Sm. Rhombicosidodecahedron  & 3.4.5.4  & $\{60,4\}$  & 207360000 & 432000 \tabularnewline
\hline 
Gt. Rhombicosidodecahedron & 4.6.10  & $\{60,3\}$  & 559872000000 & 777600000 \tabularnewline
\hline 
Snub Dodecahedron & 3.3.3.3.5  & $\{15,5\}$  & $2^{23}3^{11}5^{11}$ & $2^{20}3^{10}5^{9}$ \tabularnewline
\hline 
Trunc. Dodecahedron  & 3.10.10  & $\{30,3\}$  & 2592000 & 7200 \tabularnewline
\hline 
Gt. Rhombicuboctahedron & 4.6.4.8 & $\{24,4\}$ & 5308416 & 18432 \tabularnewline
\hline
\end{tabular}\label{T:summary}
\end{centering}
\end{table}

The minimal cover of the truncated tetrahedron is in fact $\{6,3\}_{(2,2)}$. That the latter covers the truncated tetrahedron was noted in \cite{Ha06}, but it was not shown to be a minimal cover.

\section{Analysis of presentations}

\label{sec:Analysis}

Having obtained a quotient presentation, there are a variety of questions
that one may now ask about the structure of the presentation, both
algebraically and combinatorially, that may be approached by algebraic
methods.

\subsection{Acoptic Petrie Schemes}

One such question is the determination of whether or not the given
polytope has acoptic Petrie schemes%
\footnote{Such polytopes are referred to as {\em Petrial} polytopes in \cite{Wil05}.%
}, a question related to understanding under what conditions a polyhedron
will have Petrie polygons that form simple closed curves. First, we
require some definitions; we will follow the second author's \cite{Wil05}.
A {\em Petrie polygon} of a polyhedron is a sequence of edges of
the polyhedron where any two consecutive elements of the sequence
have a vertex and face in common, but no three consecutive edges share
a common face. For the regular polyhedra, the Petrie polygons form
the equatorial skew polygons. The definition of a Petrie polygon may
be extended to polytopes of rank $n>3$ as well. An {\em exchange map} $\varrho_{i}$
is a map on the flags of the (abstract or geometric) polytope sending
each flag $\Phi$ to the unique flag 
that differs from it only by the element
at rank $i$ (this corresponds to earlier discussion of flag action for a suitable Coxeter group). A {\em Petrie map} $\sigma$ of a polytope $\mc Q$
of rank $d$ is any composition of the exchange maps $\{\varrho_{0},\varrho_{1},\ldots,\varrho_{d-1}\}$
on the flags of $\mc Q$ in which each of these maps appears exactly
once. For example, the map $\sigma=\varrho_{d-1}\varrho_{d-2}\ldots\varrho_{2}\varrho_{1}\varrho_{0}$
is a Petrie map. In particular, suppose $\mc Q\mc \simeq P(W)/N$ admits a flag action by the string C-group $W$. Then the flag action of a Coxeter element
in $W$, such as $s_{n}\ldots s_{1}s_{0}$, on a given flag in $\mc Q$ is a Petrie map.

\begin{definition} A {\em Petrie sequence} of an abstract polytope
is an infinite sequence of flags which may be written in the form
$\left(...,\,\Phi\sigma^{-1},\,\Phi,\,\Phi\sigma,\,\Phi\sigma^{2},\,...\right)$,
where $\sigma$ is a fixed Petrie map and $\Phi$ is a flag of the
polytope. \end{definition} \begin{definition} A {\em Petrie scheme}
is the shortest possible listing of the elements of a Petrie sequence. If a Petrie
sequence of an abstract polytope contains repeating cycles of elements,
then the Petrie scheme is the shortest possible cycle presentation
of that sequence. Otherwise, the Petrie scheme is the Petrie sequence.\end{definition}
For example, there is no finite presentation for a Petrie scheme of the regular tiling of the plane by squares, but while any Petrie sequence of a tetrahedron is infinitely long, any of its Petrie schemes has only four elements (and we consider cyclic permutations of a Petrie scheme to be equivalent).

A polytope possesses {\em acoptic Petrie schemes} if each proper face
appears at most once in each Petrie scheme. We borrow this terminology from Branko Gr\"unbaum who coined the term acoptic (from the Greek $\kappa o \pi\tau\omega$, to cut) to describe polyhedral surfaces with no self-intersections (cf. \cite{Gru94,Gru97,Gru99,Wil05}). Let $\{\sigma_{1},\sigma_{2},\ldots,\sigma_{t}\}$
be the collection of {\em distinct} Coxeter elements in $W$ (we assume here that $W$ is finite), and
choose $\{u_{1}=1,u_{2},u_{3},\ldots,u_{|W:N|}\}$ such that $\{\Phi^{u_{1}}=\Phi,\Phi^{u_{2}},\ldots,\Phi^{u_{|W:N|}}\}={\mc F}(\mc P(W)/N)$. Note that all Coxeter elements in $W$ are conjugates since the covering Coxeter group has a string diagram.
Following \cite{Har99}, we denote by $H_{i}$ the parabolic subgroups of $W$ of the form $\langle s_{j}:j\ne i\rangle$.
Since faces of the polytope are in one-to-one correspondence with
double cosets of the form $Nu_{j}H_{i}$, and the flag action of an
element $v\in W$ sends a face $Nu_{j}H_{i}$ in flag $\Phi^{u_{j}}$
to the face $Nu_{j}vH_{i}$ (see \cite{Har99a}), it suffices to consider
the conditions under which $Nu_{j}(\sigma_{l})^{k}H_{i}=Nu_{j}H_{i}$.
In this instance, $u_{j}(\sigma_{l})^{k}\in Nu_{j}H_{i}$, so there
exist $n\in N,h\in H_{i}$ such that $nu_{j}h=u_{j}(\sigma_{l})^{k}$.
In other words, $u_{j}^{-1}nu_{j}=(\sigma_{l})^{k}h^{-1}$, which
is equivalent to $(\sigma_{l})^{k}H_{i}\cap N^{u_{j}}\ne\emptyset$.
Note that this intersection condition depends not on our choice of
$u_{j}$, but only on the conjugates of $N$. In other words, by Theorem
\ref{th:flagCount}, we may restrict our attention only to a subcollection
of the $\Phi^{u_{j}}$, one taken from each automorphism class. Therefore,
a Petrie scheme fails to be acoptic precisely when $(\sigma_{l})^{k}H_{i}\cap N^{u_{j}}\ne\emptyset$
and $k$ is less than the size of the orbit of $\Phi^{u_{j}}$ under
the action of $\sigma_{l}$. We have thus shown the following theorem.

\begin{theorem} Let $\{u_{1}=1,u_{2},u_{3},\ldots,u_{r}\}$ be chosen
such that $\{\Phi^{u_{j}}:1\le j\le r\}$ are representatives of each
of the $r$ transitivity classes of flags under the automorphism group of the polytope $\mc P(W)/N$.
Let $\{\sigma_{1},\sigma_{2},\ldots,\sigma_{t}\}$ be the collection
of {\em distinct} Coxeter elements in $W$ and let ${m_{j,l}=|\{\Phi^{u_{j}}\alpha:\alpha\in\langle\sigma_{l}\rangle\}|}$.
Then $\mc P(W)/N$ has acoptic Petrie schemes if ${(\sigma_{l})^{k}H_{i}\cap N^{u_{j}}=\emptyset}$
for all $1\le k<m_{j,l}$. \label{th:Petrie} \end{theorem}

The results of applying such a test to the sporadic Archimedean solids are given in Table \ref{t:Petrial}.
This expands the list of known polytopes with acoptic Petrie schemes
given in \cite{Wil05} to include eight of the sporadic Archimedean
polyhedra. We say that a polytope has {\em acoptic Petrie schemes at rank $i$} if each face of rank $i$ appears at most once in each Petrie scheme, so a polyhedron has acoptic Petrie schemes if it has acoptic Petrie schemes at ranks 0, 1 and 2.%
\begin{table}[htdp]
 
\caption{The ranks at which the Archimedian polyhedra have acoptic Petrie schemes.}
\begin{centering}
\begin{tabular}{|l|c|}
\hline 
\textbf{Polyhedron} & \textbf{Acoptic Ranks}\tabularnewline
\hline
\hline 
Cuboctahedron & $\{0,1,2\}$\tabularnewline
\hline 
Great Rhombicosidodecahedron & $\{0,1,2\}$\tabularnewline
\hline 
Great Rhombicuboctahedron & $\{0,1,2\}$\tabularnewline
\hline 
Icosadodecahedron & $\{0,1,2\}$\tabularnewline
\hline 
Small Rhombicosidodecahedron & $\{0,1,2\}$\tabularnewline
\hline 
Small Rhombicuboctahedron & $\{0,1,2\}$\tabularnewline
\hline 
Pseudorhombicuboctahedron & $\emptyset$\tabularnewline
\hline 
Snub Cube & $\emptyset$\tabularnewline
\hline 
Snub Dodecahedron & $\emptyset$\tabularnewline
\hline 
Truncated Cube & $\{0,1\}$\tabularnewline
\hline 
Truncated Dodecahedron & $\{0,1\}$\tabularnewline
\hline 
Truncated Icosahedron & $\{0,1,2\}$\tabularnewline
\hline 
Truncated Octahedron & $\{0,1,2\}$\tabularnewline
\hline 
Truncated Tetrahedron & $\{0,1\}$\tabularnewline
\hline
\end{tabular}
\par\end{centering}

\label{t:Petrial} 
\end{table}

As a practical matter, one need not check all of the distinct Coxeter
elements, but instead only half of them, since the inverse of a Coxeter
element is itself a Coxeter element, and inverse pairs generate the
same sequences of flags, only in reverse order. Thus for polyhedra,
one need only check $\sigma_{1}=s_{0}s_{1}s_{2}$ and $\sigma_{2}=s_{0}s_{2}s_{1}$.

Let $|\sigma_{l}|$ denote the order of $\sigma_{l}$. It is worth noting that it is easy to construct examples of polytopes
for which $m_{j,l}<|\sigma_{l}|$ for all $j$ and $l$, even when
the covering regular polytope is finite and all of the schemes are
acoptic. One such is obtained by taking the quotient of the universal square tessellation
$\{4,4\}$, whose automorphism group $W$ is the Coxeter group $[4,4]$. Now let $N=\langle(\nu_{1}\nu_{2})^{3},(\nu_{1}\nu_{2}^{-1})^{5}\rangle$
where $\nu_{1}=s_{0}s_{1}s_{2}s_{1}$ and $\nu_{2}=s_{1}s_{0}s_{1}s_{2}$. Then $\mc P(W)/N=[4,4]/N$ is a toroidal polyhedron. In this case, $m_{j,l}$ is
either 6 or 10, but $|\sigma_{l}|=30$ in $W/Core(W,N)$. For a further
discussion of Petrie polygons and polytopes with acoptic Petrie schemes
see \cite{Wil05}.

\subsection{Size of Presentations}

The pseudorhombicuboctahedron (also known as the elongated square
gyrobicupola, or Johnson solid $J_{37}$)\footnote{The pseudorhombicuboctahedron has been ``discovered'' independently on numerous occasions and has proved to be an excellent example of the difficulties mathematicians have in constructing definitions about intuitively understood objects that are sufficiently rigorous so as to specify precisely the objects they wish to study without accidentally assuming unstated constraints (such as symmetry). The interested reader is encouraged to review Gr\"unbaum's excellent discussion of the history in \cite{Gru08}.}
  provides an interesting
case for discussion, because while it has the same local structure
as the small rhombicuboctahedron (vertex stars of type $3.4.4.4$),
it has significantly less symmetry. Theorem \ref{th:flagCount} provides
a computationally very fast method of determining that there are in
fact twelve equivalence classes of flags (a fact otherwise tedious
to determine), while Theorem \ref{th:Petrie} provides a rapid method
of verifying that the Petrie schemes of $J_{37}$ are not all acoptic
at any rank. Perhaps more surprising to the reader might be the comparison
of the sizes of the group presentation with the small rhombicuboctahedron.
While the minimal cover of the small rhombicuboctahedron is of order
$1\,327\,104$ the cover for the pseudorhombicuboctahedron is more than
ten orders of magnitude larger at $16\,072\,626\,615\,091\,200$.

\section{Some Open Questions}

\label{sec:Conclusion} We include here some questions motivated by
the current work. Theorem \ref{th:core} provides a minimal presentation
for a polytope as a quotient of a regular polytope, but only in the
instance that $\mc P(W/\Core(W,N))$ is a well defined polytope. Does
there exist an example of a (finite) polytope for which $\mc P(W/\Core(W,N))$
is not polytopal? Also, in the examples studied to date, finite polytopes
have all yielded representations as the quotients of finite regular
polytopes. Is there an example of a finite polytope which does not
admit a presentation as the quotient of a finite regular polytope?
Both of these questions would be answered in the negative if the following
conjecture --- and thus its corollary by Theorem \ref{th:flagCount}
--- are true (for definitions and a more detailed discussion of the
role semisparse subgroups play in the theory of quotient representations,
see \cite{Ha06}). 

\begin{conjecture}
If $N$ is semisparse in $W$ then $\Core(W,N)$ is also semisparse.\label{th:semiSparse}
\end{conjecture}
\begin{corollary} Assuming Conjecture \ref{th:semiSparse}, every finite abstract polytope admits a presentation
as the quotient of a finite regular abstract polytope.\end{corollary}
A computer survey of the symmetry groups of abstract regular polytopes
found no counterexamples to Conjecture \ref{th:semiSparse} for groups
$W$ of order less than 639.

\bibliographystyle{amsalpha} 
\bibliography{ArchimedeanBibliography}

\end{document}